\documentclass[11pt,reqno]{amsart}

\renewcommand{\P}{\mathbb{P}}
\newcommand{\abs}[1]{\left\vert #1 \right\vert}
\newcommand{\diff}{\mathop{}\!\mathrm{d}}

\usepackage{epsfig,amssymb,amsmath,amscd,amsfonts,color, amsthm}
\title[Large deviations for the largest eigenvalue]{Large Deviations
  for the largest eigenvalue of rank one deformations of Gaussian ensembles}

\author[M. Maida]{Myl\`ene Ma\"{\i}da}
\thanks{Universit\'e Paris-Sud
Laboratoire de math\'ematiques
 91405 Orsay Cedex, France.\\
\indent Email: mylene.maida@math.u-psud.fr}

\def\H{H}
\newtheorem{prop}{Proposition}[section]
\newtheorem{theo}[prop]{Theorem}
\newtheorem{defi}[prop]{Definition}
\newtheorem{lem}[prop]{Lemma}
\newtheorem{fact}[prop]{Fact}

\newtheorem{cor}[prop]{Corollary}

\newtheorem{rmk}[prop]{Remark}

\def\P{\Bbb P}
\def\E{\Bbb E}

\newcommand{\RR}{\mathbb{R}}

\newcommand{\NN}{\mathbb{N}}
\newcommand{\PP}{\mathbb{P}}
\newcommand{\EE}{\mathbb{E}}

\newcommand{\VV}{\mathbb{V}}
\newcommand{\PR}{\mathcal{P}(\RR)}

\newcommand{\ppq}{\leqslant}
\newcommand{\pgq}{\geqslant}
\renewcommand{\ge}{\geqslant}
\renewcommand{\geq}{\geqslant}
\renewcommand{\le}{\leqslant}
\renewcommand{\leq}{\leqslant}

\newcommand{\muN}{\hat{\mu}_N}

\newcommand{\pin}{\hat{\pi}_N}

\def\tr{{\mbox{tr}}}

\def\QNA{\mathbb Q_N^\theta}
\def\PNA{\mathbb P_N^\theta}
\def\ZNA{Z_N^\theta}

\def\b{\beta}
\def\d{\delta}

\def\k{\kappa}
\def\l{\lambda}

\def\s{\sigma}
\def\t{\theta}

\def\D{\Delta}

\def\ra{\rightarrow}

\def\Aa{{\mathcal A}}

\def\Ca{{\mathcal C}}

\def\Ja{{\mathcal J}}

\def\Oa{{\mathcal O}}

\def\Ua{{\mathcal U}}

\def\nun{{\hat\nu_N}}

\def\Pa{{\mathcal{P}}}

\def\H{H_{\mu}}

\def\R{R_{\mu}}

\renewcommand{\epsilon}{\varepsilon}
\renewcommand{\phi}{\varphi}

\newlength{\aux}
\def\KBT{K_\theta^\beta}
\def\FBT{F_\theta^\beta}

\begin{document}

\begin{abstract}
We establish a large deviation principle for the largest eigenvalue of
a rank one deformation of a matrix from the GUE or GOE. As a
corollary, we get another proof of the phenomenon, well-known in 
learning theory and finance, that the largest eigenvalue separates
from the bulk when the perturbation is large enough.\\
A large part of the paper is devoted to an auxiliary result
on the continuity of spherical integrals in the case when one of the
matrix  is of rank one, as studied in \cite{GM2}.

\end{abstract}

\maketitle


\section{Introduction}

We consider in this paper 
rank one deformations of Wigner matrices, 
that is matrices which can be written 
 $W_N+A_N,$ with $W_N$
from the Gaussian Orthogonal (or Unitary) Ensemble
and $A_N$ rank one deterministic.

Since the fifties, the classical Gaussian ensembles (see Mehta \cite{Meh})
have been extensively studied. Various results for the global regime
were established (Wigner semicircle law \cite{Wig}, large deviations for the
spectral measure \cite{BA-G}...); the statistics of the spacings between 
eigenvalues were investigated \cite{Dei, DKM}, as well as the behaviour of extremal
eigenvalues
(Tracy-Widom distribution \cite{TW}). In the meantime, people got
interested in the 
universality of some of these results. In this context, it is natural to
look at various deformations of these ensembles, for example
the rank one deformations we are interested in.

This so-called ``deformed Wigner ensemble'' was studied
in  \cite{Joh} and  \cite{BH},  where the authors focused mainly on the problem of the
local spacings and in \cite{SP} and \cite{FP}, where they studied
 the behaviour of the largest eigenvalue. In this framework,
our goal in this paper will be to establish a large deviation
principle
for the largest eigenvalue of $X_N= W_N+A_N$, that we denote in the
sequel by $x_N^*$. Note that our result can also be seen as a generalization of
the result established in \cite{BDG} for the largest eigenvalue of a
matrix
distributed according to the GOE. \\
If we denote by $\t$ the unique non zero eigenvalue of $A_N,$
the joint law of the eigenvalues $x_1, \ldots, x_N$
of $X_N = W_N+A_N$ is given by 
\begin{equation}
\label{QNA}
\QNA(dx_1, \ldots, dx_N) = \frac 1 \ZNA 
\prod_{i <j} |x_i - x_j|^\beta I_N^\beta(\t,X_N) e^{-\frac N 2 \sum_{i=1}^N x_i^2}
dx_1 \ldots dx_N,
\end{equation}
where $I_N^\beta$ is the spherical integral defined by
$$ I_N^\beta(\t,X_N) := \int e^{N \tr (UX_NU^*A_N)} dm_N^\beta(U) =  \int e^{N \t (UX_NU^*)_{11}} dm_N^\beta(U),$$
with $m_N^\beta$ the Haar probability measure on $\Oa_N$ the 
orthogonal group of size $N$ if $\beta=1$, on the unitary group 
$\Ua_N$ if $\beta=2$ and $\ZNA$ is a normalizing constant.
The fact that the joint law of the eigenvalues of $X_N $ and that 
$I_N^\beta(\t,X_N) $ depend on $A_N$ only through 
its non zero eigenvalue $\t$ come from the unitary invariance
respectively of the law of $W_N$ and of  the Haar measure  $m_N^\beta$ .\\
 Our main result can be roughly stated as follows
\begin{theo}
\label{pgdflou}
If $\t \pgq 0,$
then under $\QNA$, the largest eigenvalue $x_N^*= \max\{x_1, \ldots, x_N\} $
satisfies a large deviation principle in the scale $N,$
with good rate function $K_\theta^\beta.$
\end{theo}

We get an explicit expression of $K_\theta^\beta$ on which   one can
see in particular that it differs from the rate function for the
deviations of the largest eigenvalue of the non-deformed
model that was obtained in \cite{BDG}.\\

Note that in the case when $\t <0,$ similar results
would hold  for the smallest eigenvalue of the deformed ensemble. 
We let the precise statement to the reader and assume in the sequel
that
$\t >0.$

We have to mention 
an important corollary  of Theorem \ref{pgdflou} :
\begin{cor}
\label{seppgv}
Under $\QNA$, $x_N^*$ converges almost surely to
 the edge of the support of the semicircle law $\s_\b$ as long as
$\t \ppq \t_c := \sqrt{\frac{\beta}{2}}$ and  separates from the
 support
 when $\t > \t_c.$ In this case, it converges 
to $\t + \frac{\beta}{2\t}.$ 
\end{cor}

This allows us to give a new proof, via large deviations, to this
known phenomena which is crucial for applications to finance and
learning theory (cf. for example \cite{HR, LCPB}).\\

The organisation of the paper is as follows : 
as we can see in \eqref{QNA} above, the expression of the joint law
$\QNA$ of the eigenvalues involves the spherical integrals
$ I_N^\beta$ in the case when one of the matrices is of rank one.
We got the asymptotics
of this quantity in  \cite{GM2} 
but we will need a precise continuity result of these
spherical integrals to whom Section \ref{sec:cont} is devoted.
In Section \ref{sec:ldp}, we give in Theorem \ref{main} a precise statement of Theorem
\ref{pgdflou} and prove it. Finally, in a very short Section \ref{sec:sep},
we show how to derive Corollary \ref{seppgv} from this Large 
Deviation Principle.

\section{Continuity of spherical integrals}
\label{sec:cont}
The question we want to address in this section is the continuity,
in a topology to be prescribed, of $ I_N^\beta(\t,B_N),$ in its
second argument $B_N.$

We denote by $\lambda_1(B_N), \ldots,\lambda_N(B_N) $ the eigenvalues of $B_N$
in decreasing order and we let \linebreak
 $\displaystyle \hat \nu_{B_N} := \frac 1 N
\sum_{i=2}^N \d_{\lambda_i(B_N)}$ and $\|B_N\|_\infty:= \max\left(|\lambda_1(B_N)|,|\lambda_N(B_N)|\right)$; $d$ is the Dudley distance defined on
probability measures by 
$$ d(\mu,\nu) = \sup \left\{\left|\int fd\mu - \int f d\nu\right| ; |f(x)| \ppq 1 \textrm{ and  }
\left| \frac{f(x)-f(y)}{x-y}\right| \ppq 1, \forall x\neq y\right\}.$$
The following continuity property holds

\begin{prop}
\label{cont}
For any $\t>0,$  for any 
$\k \in ]0, 1/2[,$ there exists a function $g_\k: \RR_+ \ra \RR_+$ going
    to zero at zero such that, for any $\d>0 $ and $N$ large enough,
if $B_N$ and $B_N^\prime$ are such that
$d(\hat \nu_{B_N},\hat \nu_{B_N^\prime})\ppq N^{-\k}$ and
$|\lambda_1(B_N)-\lambda_1(B_N^\prime)|\ppq \d,$ with $\sup_N \|B_N\|_\infty
<\infty$ then

$$\left| \frac 1 N \log  I_N^\beta(\t,B_N) - \frac 1 N \log
I_N^\beta(\t,B_N^\prime) \right| \ppq g_\k(\d).$$
\end{prop}

\begin{rmk}
According to Theorem 6 of \cite{GM2}, we know that, for some values of
$\t$,
the limit of $\frac 1 N \log  I_N^\beta(\t,B_N)   $
as $N$ goes to infinity depends not only on the limiting spectral
measure of $B_N$ but also on the limit of $\lambda_1(B_N).$
Therefore  $\frac 1 N \log  I_N^\beta(\t,B_N)   $ cannot be
continuous in the spectral
measure of $B_N$ but we have also to localize  $\lambda_1(B_N).$
That is precisely the content of Proposition \ref{cont} above.
We also refer the reader to the remarks made in \cite{GM2} on
point $(3)$ of Lemma 14 therein.
\end{rmk}

A key step to show  Proposition \ref{cont} is to get an equivalent as
explicit as possible of $\frac 1 N \log  I_N^\beta(\t,B_N)  . $ This
is given by
\begin{lem}\label{approx}
Fix $\theta > 0$, and suppose that $(B_N)_{N\ge 1}$ is a sequence of diagonal matrices such that 
\[
    \|B\|_\infty = \sup_N   \|B_N\|_\infty< \infty.
\]
Then for every $\epsilon > 0$ and for every $N \geq N_0(\theta,\|B\|_\infty,\epsilon)$ we have 
\[
    \left|\frac{1}{N} \log I_N(\theta,E_N) - \left(\theta v_N - \frac{\beta}{2N} \sum_{i=1}^N \log \left(1+2\theta v_N - 2\theta \lambda_i \right)\right)\right| \leq C(\theta, \|B\|_\infty)N^{-\frac{1}{3}+\epsilon}
\]
where $v_N$ is the unique solution in $[\lambda_N(B_N)-\frac{\beta}{2\theta}, \lambda_1(B_N) - \frac{\beta}{2\theta}]^c$ to the equation
\[
    \frac{1}{N} \frac{\beta}{2\theta} \sum_{i=1}^N \frac{1}{v_N + \frac{\beta}{2\theta} - \lambda_i} = 1.
\]
\end{lem}

This lemma can be regarded as a generalization to any value of $\t$ of
the second point of Lemma 14 in \cite{GM2}. \\
The remaining of this section is devoted to its proof. 
For sake of simplicity, we prove in full details the case $\b =1$ and
let
to the reader the changes to the other cases.\\

\subsection{Some preliminary results}
\mbox{}\\

We first introduce some notations.\\

\textbf{Notations.}
\begin{itemize}
\item We denote by  $\lambda_1\pgq \ldots \pgq \lambda_N$  the eigenvalues of
  $B_N$
in decreasing order.
\item $\EE$ and $\VV$ denotes respectively the expectation and the variance
  under the standard Gaussian measure $\PP$ on $\RR^N.$
\item As    $1+2\t v_N -2\t \lambda_1 >0$, we can define a probability
  measure given by 
$$ P_N(dg_1,\ldots, dg_N) = (2\pi)^{-\frac{N}{2}} \prod_{i=1}^N
\left[ \sqrt{1+2\t v_N-2\t \lambda_i} \,\,e^{-\frac 1 2 (1+2\t v_N-2\t \lambda_i)g_i^2}
  dg_i\right].$$
We denote by $\EE_{P_N}$ and $\VV_{P_N}$ respectively the
expectation and the variance under $P_N.$
\end{itemize}

Before going to the proof of Lemma \ref{approx}, we enumerate
hereafter
some inequalities on the quantities we have just introduced, that will
be useful further.

\begin{fact} \label{ineq}
We have the following inequalities :
\[ v_N \le \lambda_1 \quad \textrm{ and } \quad \forall 1 \le i \le N,   v_N +\frac 1{2\t}
  -\lambda_i \pgq \frac 1 {2\t N}.\]

%
\end{fact}

\begin{proof}

The function $x \mapsto \frac{1}{N}\sum_{i=1}^N \frac 1 { x +\frac 1{2\t} -\lambda_i} - 2\theta$ is decreasing on $(\lambda_1 - \frac{1}{2\t}, \infty),$ converges to $\infty$ as $x$ decreases to $\lambda_1 - \frac{1}{2\t}$ and is nonpositive at $x=\lambda_1.$  
Therefore, it vanishes at $v_N \le \lambda_1.$

As for the second inequality, we have that
$ \sum_{i=1}^N \frac 1 { v_N +\frac 1{2\t} -\lambda_i} = 2 \t N$ and each
term is positive so that any of them  is smaller than $2\t N.$ 
\end{proof}

We also need the following auxiliary result :
\begin{lem}
\label{lem:chi_concentration}
Let $\{X_i\}_{i=1}^N$ be independent standard Gaussian random variables  under the measure $\P$ on $\mathbb R^N$, and suppose that $\alpha_1, \ldots, \alpha_N$ are positive numbers summing to one. Then for any $0 < \kappa < \frac{1}{2}$ there exists $C(\kappa) > 0$ such that, for all $\epsilon > 0$ sufficiently small and for all $N \geq N_0(\kappa,\epsilon)$
\[
    \frac{1}{N}\log\P \left( \abs{\sum_{i=1}^N \alpha_i (X_i^2-1)} \leq N^{-\kappa} \right) \geq -C(\kappa)\frac{N^{2\kappa+\epsilon}}{N}
\]
\end{lem}
\begin{proof}
Let $\epsilon_N \to 0$ to be chosen. We have
\begin{align*}
   \P\left( \abs{ \sum_{k=1}^N \alpha_k(X_k^2-1)} \leq N^{-\kappa} \right) &\geq\P\left( \abs{ \sum_{k : \alpha_k < \epsilon_N} \alpha_k(X_k^2 - 1)} \leq \frac{N^{-\kappa}}{2}, \abs{ \sum_{k : \alpha_k \geq \epsilon_N} \alpha_k(X_k^2 - 1)} \leq \frac{N^{-\kappa}}{2} \right) \\
    &=\P\left( \abs{ \sum_{k : \alpha_k < \epsilon_N} \alpha_k(X_k^2 - 1)} \leq \frac{N^{-\kappa}}{2} \right)\P \left(\abs{ \sum_{k : \alpha_k \geq \epsilon_N} \alpha_k(X_k^2 - 1)} \leq \frac{N^{-\kappa}}{2} \right)
\end{align*}
We first bound the contribution from the small $\alpha_k$ values: Since
\[
    \mathrm{Var}\left( \sum_{k : \alpha_k < \epsilon_N} \alpha_k(X_k^2-1) \right) = 2\sum_{k : \alpha_k < \epsilon_N} \alpha_k^2 \leq 2 \epsilon_N \sum_{k : \alpha_k < \epsilon_N} \alpha_k \leq 2\epsilon_N,
\]
Chebyshev's gives us
\[
   \P\left( \abs{ \sum_{k : \alpha_k < \epsilon_N} \alpha_k(X_k^2 - 1)} > \frac{N^{-\kappa}}{2} \right) \leq 8N^{2\kappa}\epsilon_N.
\]
If $\epsilon_N = N^{-2\kappa-\epsilon}$ for some small $\epsilon > 0$, for example, then this is $o(1)$; hence for $N$ sufficiently large depending on $\kappa$ we have
\[
   \P\left( \abs{ \sum_{k : \alpha_k < \epsilon_N} \alpha_k(X_k^2 - 1)} \leq \frac{N^{-\kappa}}{2} \right) \geq \frac{1}{2}.
\]

We now bound the contribution from the large $\alpha_k$'s. Let $m_N$ be the number of $\alpha_k$'s greater than $\epsilon_N$. Then $1 = \sum_{k} \alpha_k \geq \sum_{k : \alpha_k \geq \epsilon_N} \alpha_k \geq m_N\epsilon_N$, so $m_N \leq \frac{1}{\epsilon_N}$ and thus for some universal constant $c$ we have
\begin{align*}
   \P \left(\abs{ \sum_{k : \alpha_k \geq \epsilon_N} \alpha_k(X_k^2 - 1)} \leq \frac{N^{-\kappa}}{2} \right) &\geq\P \left(\forall \, \alpha_k \geq \epsilon_N, \abs{\alpha_k(X_k^2-1)} \leq \frac{N^{-\kappa}}{2m_N} \right) \\
    &= \prod_{k : \alpha_k \geq \epsilon_N}\P\left(\abs{X_k^2-1} \leq \frac{N^{-\kappa}}{2m_N\alpha_k} \right) \\
    &\geq \prod_{k : \alpha_k \geq \epsilon_N}\P \left( \abs{X_k^2-1} \leq \frac{N^{-\kappa}}{2m_N} \right) \\
    &\geq \left(\P \left( \abs{X_k^2-1} \leq \frac{N^{-\kappa}\epsilon_N}{2} \right) \right)^{\frac{1}{\epsilon_N}} \\
    &\geq \left(c\frac{N^{-3\kappa-\epsilon}}{2}\right)^{N^{2\kappa+\epsilon}}
\end{align*}
Hence
\begin{align*}
    \frac{1}{N}\log\P \left( \abs{\sum_{i=1}^N \alpha_i (X_i^2-1)} \leq N^{-\kappa} \right) &\geq -\frac{\log2}{N} +
    \frac{1}{N}\log\left(\P \left(\abs{ \sum_{k : \alpha_k \geq \epsilon_N} \alpha_k(X_k^2 - 1)} \leq \frac{N^{-\kappa}}{2} \right) \right) \\
    &\geq -\frac{\log 2}{N} + \frac{N^{2\kappa+\epsilon}}{N} \log\left( \frac{c}{2}N^{-3\kappa-\epsilon} \right) \\
    & \geq  -\frac{\log 2}{N} -(3\kappa+2\epsilon) \frac{N^{2\kappa+\epsilon}\log(N)}{N},
\end{align*}
for $N$ large enough, which gives the result by choosing $\epsilon$ such that $2\kappa+\epsilon < 1$. 

\end{proof}

\subsection{Proof of Lemma \ref{approx}}
\mbox{}\\
We can now go back to the \textbf{proof of the upper bound:} 
the starting point will be the same as in \cite{GM2}. It is a well
known
fact that the first column vector of a random orthogonal matrix
distributed according to the Haar measure on $\Oa_N$ has the same
law as a standard Gaussian vector in $\RR^N$ divided by its Euclidian
norm.
Therefore, we
can write
$$ I_N(\t, B_N) = \EE \left( \exp\left\{ N\t \frac{\sum_{i=1}^N \lambda_i
  g_i^2}{\sum_{i=1}^N g_i^2}\right\}\right).$$
From concentration for the norm of a Gaussian vector (cf.  \cite{GM2}
for details), we get that, for any $\k$ such that $0<\k<1/2,$
\begin{equation} \label{encadr}
1\ppq \frac{ I_N(\t, B_N)}{\EE \left( \mathbf{1}_{\Aa_N(\k)}\exp\left\{ N\t \frac{\sum_{i=1}^N \lambda_i
  g_i^2}{\sum_{i=1}^N g_i^2}\right\}\right)}\ppq \d(\k,N),
\end{equation}
where $\Aa_N(\k) = \left\{\left|\frac{\|g\|^2}N -1\right|\ppq N^{-\k}\right\}$
and $ \d(\k,N)$ goes to one at infinity for any  $0<\k<1/2.$\\

From there, we have 
\begin{eqnarray*}
 I_N(\t, B_N) & \ppq & \d(\k,N) e^{N\t v_N+ N^{1-\k}\t(\|B\|_\infty+v_N)}
\EE\left[ \mathbf 1_{\Aa_N(\k)} \exp\left\{ \t \sum_{i=1}^N \lambda_i g_i^2
- \t v_N   \sum_{i=1}^N g_i^2 \right\}\right]\\
& \ppq & \d(\k,N) e^{N\t v_N+ N^{1-\k}\t (\|B\|_\infty+v_N)} \prod_{i=1}^N
\left[ \sqrt{1+2\t v_N-2\t \lambda_i}\right]^{-1}.
\end{eqnarray*}
 Therefore, for any $0<\kappa <1/2,$
we get that
for $N$
large enough,
$$
\frac 1 N \log  I_N(\t, B_N) - \left(\t v_N - \frac 1 2 
 \sum_{i=1}^N \log\left( 1+2\t v_N-2\t \lambda_i\right)\right) \le \frac{\t (\|B\|_\infty+v_N)}{N^\kappa} +\frac{\delta(\kappa,N)}{N}.$$
Now Fact 2.4 tells us that $v_N \leq \lambda_1 \leq \|B\|_\infty$, and $v_N \geq \lambda_1 - \frac{1}{2\theta} \geq -(\|B\|_\infty+\frac{1}{2\theta})$, so $\abs{v_N} \leq \|B\|_\infty+\frac{1}{2\theta}$. Hence
\[
    \frac{1}{N} \log I_N(\theta,E_N) - \left(\theta v_N - \frac{1}{2N} \sum_{i=1}^N \log \left(1+2\theta v_N - 2\theta \lambda_i \right)\right) \leq \frac{\theta(2\|B\|_\infty + \frac{1}{2\theta})}{N^{\kappa}} + \frac{\log 2}{N} \leq C(\theta, \|B\|_\infty) N^{-\kappa}
\]
where the last inequality holds for $N$ sufficiently large. 

 For the lower bound, we start with the following bound that can be deduced from \eqref{encadr}:
\begin{align*}
    \frac{1}{N}I_N(\theta,E_N) - \left(\theta v_N - \frac{1}{2N} \sum_{i=1}^N \log \left(1+2\theta v_N - 2\theta \lambda_i \right)\right) &\geq -\frac{\theta(\abs{v_N}+\|B\|_\infty}{N^{\kappa}} + \frac{1}{N}\log P_N(A_N(\kappa)) \\
    &\geq -\frac{\theta(2\|B\|_\infty+ \frac{1}{2\theta})}{N^{\kappa}} + \frac{1}{N}\log P_N(A_N(\kappa)).
\end{align*}
Here, $P_N$ is the probability measure on $\mathbb R^N$ given by
\[
    P_N(\diff g_1, \ldots, \diff g_N) = \frac{1}{\sqrt{2\pi}^N} \prod_{i=1}^N \left[ \sqrt{1+2\theta v_N - 2\theta \lambda_i} e^{-\frac{1}{2}(1+2\theta v_N - 2\theta \lambda_i) g_i^2} \diff g_i \right],
\]
and $A_N(\kappa)$ is the event $\{\abs{\frac{\|g\|^2}{N} - 1} \leq N^{-\kappa}\}$. 

Set $\alpha_i = \frac{1}{N}\frac{1}{1+2\theta v_N - 2\theta \lambda_i}$. By the definition of $v_N$, each $\alpha_i$ is positive and $\sum_{i=1}^N \alpha_i = 1$. Thus an application of Lemma \ref{lem:chi_concentration} gives us
\[
    \frac{1}{N}\log P_N(A_N(\kappa)) = \frac{1}{N}\log \P_{\text{i.i.d.}} \left( \abs{\sum_{i=1}^N \alpha_i(X_i^2 - 1) } \leq N^{-\kappa} \right) \geq -C(\kappa)\frac{N^{2\kappa+\epsilon}}{N}.
\]
Thus
\[
    \frac{1}{N}I_N(\theta,E_N) - \left(\theta v_N - \frac{1}{2N} \sum_{i=1}^N \log \left(1+2\theta v_N - 2\theta \lambda_i \right)\right) \geq -C(\theta,\|B\|_\infty N^{-\kappa} - C(\kappa)N^{2\kappa+\epsilon-1}.
\]
These rates of decay compete; the optimal rate of decay, which is $N^{-\frac{1}{3}+\epsilon}$, occurs when we choose $\kappa = \frac{1}{3}$.

This concludes the proof of Lemma \ref{approx}. \hfill $\Box$

\subsection{Proof of Proposition \ref{cont}}
\mbox{}\\

Let $\k \in \left( 0, 1/2\right)$ be fixed and $(B_N)_{N \in \NN}$
and  $(B_N^\prime)_{N \in \NN}$ two sequences of matrices. 
We denote by $\lambda_1 \pgq \ldots \pgq \lambda_N$ and $\lambda_1^\prime \pgq \ldots
\pgq \l^\prime_N$
the eigenvalues of $B_N$ and $B_N^\prime$ respectively, both in decreasing
order.\\
We assume that $\displaystyle d(\hat \nu_{B_N}, \hat \nu_{B_N^\prime}) \ppq N^{-\k}$
and $|\lambda_1 - \lambda_1^\prime | \ppq \d$ for $N$ large enough and that 
there exists $M$ such that $\sup_N \|B_N\|_\infty < M$ and 
$\sup_N \|B_N^\prime\|_\infty < M.$\\
We introduce also the following notation: 
$ \displaystyle H_{B_N}(z) = \frac 1 N \sum_{i=1}^N \frac 1 {z-\lambda_i}$
and  $\displaystyle H_{B_N^\prime}(z) = \frac 1 N \sum_{i=1}^N \frac 1
{z-\lambda_i^\prime}.$\\

\noindent
$\bullet$ \textbf{First case :} $N$ and $\t$ are such that 
$2\t \in H_{B_N}( (\lambda_1+2\d, +\infty)) \cap  H_{B_N^\prime}(
  (\lambda_1^\prime+2\d, +\infty)). $\\

\noindent
In this frame work, continuity has been established in Lemma 14
of \cite{GM2}. It comes from Lemma \ref{approx} and the fact that, for 
any $\l \in \cup_{N_0 \pgq 0} \cap_{N\pgq N_0} (\textrm{supp }\hat \nu_{B_N}
\cap \textrm{supp } \hat \nu_{B_N^\prime} ),$ $z \mapsto (z-\l)^{-1}$ is
continuous bounded (with a norm independent of $\d$) on 
$ \cup_{N_0 \pgq 0} \cap_{N\pgq N_0} ((\lambda_1+2\d, +\infty) \cap
(\l^\prime_1+2\d, +\infty)).$\\

\noindent
$\bullet$ \textbf{Second case :}  $N$ and $\t$ are such that 
$2\t \notin H_{B_N}((\lambda_1 +2\d, + \infty))$ and
$2\t \notin H_{B_N^\prime}((\l^\prime_1 +2\d, + \infty)) .$\\

\noindent
In this case, $2\t \notin H_{B_N}((\lambda_1 +2\d, + \infty))
\Rightarrow v_N +\frac{1}{2\t} \in (\lambda_1, \lambda_1 + 2\d)$ and similarly
$ v_N^\prime +\frac{1}{2\t} \in (\lambda_1^\prime, \lambda_1^\prime + 2\d)$ so that 
$$ |v_N - v_N^\prime | \ppq 3 \d.$$
Thanks to Lemma \ref{approx}, we know that it is enough  to study
$$ \D_N := 
\left| \frac{1}{N} \sum_{i=1}^N \log\left(v_N + \frac{1}{2\t}- \lambda_i\right)
- \frac{1}{N} \sum_{i=1}^N \log\left(v^\prime_N + \frac{1}{2\t}- 
\l^\prime_i\right) \right|.$$

As  $d(\hat \nu_{B_N}, \hat \nu_{B_N^\prime}) \ppq N^{-\k}$, we
proceed
as in the proof of Lemma 5.1 in \cite{GZ} 
and define a permutation $\s_N$ that allows to put in pairs
all but $(N^{1-\k}\wedge N\d)$ of the $\lambda_i$'s
with a corresponding $\lambda_{\s_N(i)}^\prime$ 
which lies at a distance less than $\d$ from $\lambda_i$.\\
As in \cite{GZ}, we denote by $\Ja_0$ the set of indices $i$
such that we have such a pairing.
Then we have 
\begin{eqnarray*}
\D_N & \ppq & \frac 1 N \sum_{i \in \Ja_0}
\max\left(\frac{1}{v_N + \frac{1}{2\t}- \lambda_i}, 
\frac{1}{v^\prime_N + \frac{1}{2\t}- \lambda^\prime_{\s_N(i)}}\right)
(|v_N - v^\prime_N|+|\lambda_i - \lambda^\prime_{\s_N(i)}|) \\
& & +  \frac 1 N \sum_{i\in \Ja_0^c} \left| \log\left(v_N + \frac{1}{2\t}- \lambda_i\right)
- \log\left(v^\prime_N + \frac{1}{2\t}- 
\lambda^\prime_i\right)\right|\\
& \ppq & \left(\frac 1 N \sum_{i=1}^N\frac{1}{v_N + \frac{1}{2\t}- \lambda_i} +
\frac 1 N \sum_{i=1}^N \frac{1}{v^\prime_N + \frac{1}{2\t}- \lambda^\prime_{\s_N(i)}}\right) 
 4\d \\
& & +  \frac 1 N \sum_{i\in \Ja_0^c} \left| \log\left(v_N + \frac{1}{2\t}- \lambda_i\right)
-  \log\left(v^\prime_N + \frac{1}{2\t}- 
\lambda^\prime_i\right)\right|\\
& \ppq & 16 \t\d + \frac 2 N [N^{1-\k} \wedge N\d] |\log(2N \t)|\vee
\log\left(2M+\frac 1 {2\t}\right),
\end{eqnarray*}
where we used once again that 
$$ \frac 1 {2N\t} \ppq v_N + \frac{1}{2\t}- \lambda_i \ppq  2M+\frac 1 {2\t}$$
so that we get the required continuity in this second case.\\

\noindent
$\bullet$ \textbf{Third case :}  $N$ and $\t$ are such that 
$2\t \in H_{B_N}((\lambda_1 +2\d, + \infty))$ and
$2\t \notin H_{B_N^\prime}((\lambda^\prime_1 +2\d, + \infty)) .$\\

In this case, we proceed exactly as in the second case. The only point
is that establishing that $v_N$ cannot be far from $v_N^\prime$ will
be a bit more involved. We address this point in detail.\\

On one side we have from Fact \ref{ineq} that
\begin{equation}
\label{vun}
 \l^\prime_1 \ppq v^\prime_N  + \frac{1}{2\t}\ppq \l^\prime_1 + 2\d.
\end{equation}
On the other side, as  $|\lambda_1 - \lambda_1^\prime| \ppq \d$, $\lambda_1^\prime + 2\d$ 
is greater than $\lambda_1$ and
the map $B_N \mapsto H_{B_N}$ is continuous outside 
the support of all the spectral measures
so that
$$ \left| H_{B_N^\prime}(\l^\prime_1 + 2\d) - 
H_{B_N}(\l^\prime_1 + 2\d)\right|\ppq C(N,\d),$$
with, for any fixed $\delta >0,$ the function $C(N,\delta)$ going to zero as $N$ goes to infinity.\\

Furthermore, $H_{B_N^\prime}$ is decreasing on $(\l^\prime_1, + \infty)$
so that $H_{B_N^\prime}(\l^\prime_1 + 2\d) < 2\t$, yielding 
$$ H_{B_N}(\l^\prime_1 + 2\d) \ppq 2\t + C(N,\d),$$
and $H_{B_N}$ being decreasing
$$H_{B_N}(\lambda_1 + 3\d) \ppq 2\t + C(N,\d) = H_{B_N} \left(v_N + \frac{1}{2\t}\right) +  C(N,\d),$$
what implies 
$$ \lambda_1 \ppq v_N +\frac{1}{2\t} \ppq \lambda_1 + 3\d + K(N,\d), $$
with, for any fixed $\delta >0,$ the function $K(N,\delta)$ going to zero as $N$ goes to infinity.
and, together with \eqref{vun} this gives that
$$ |v_N - v_N^\prime | \ppq 5 \d + K(N,\d).$$

Now  the same estimates as in the second case above 
lead to the same conclusion.
This gives Proposition \ref{cont}. \hfill $\Box$

\section{Large deviations for $x_N^*$}
\label{sec:ldp}

The goal of this section is to give a precise statement and then prove
a large deviation principle for  $x_N^*$, the largest eigenvalue of a
matrix from the deformed Gaussian ensemble. In order to do that, we have
to recall a few notations.

\begin{defi}
For $\mu$ a compactly supported measure, we define $H_\mu$ 
its Hilbert transform by
\begin{eqnarray*}
H_\mu : \RR \setminus co(\textrm{supp } \mu) & \rightarrow  &\RR\\
z & \mapsto & \int \frac 1 {z-\l} d\mu(\l).
\end{eqnarray*}
 with $co(\textrm{supp }  \mu)$ the convex enveloppe 
of the support of $\mu$.\\
It is easy to check that $H_\mu $ is injective, therefore
we can define its functional inverse $G_\mu,$ and the $R$-transform
$R_\mu$ is given, for $z \neq 0,$ 
by $R_\mu(z) = G_\mu(z) -\frac 1 z$. Note that it  can be analytically continued
at $0$.
\end{defi}

\noindent
We can now state

\begin{theo} \label{main}
Under the measure 
$$ \QNA(dx_1, \ldots, dx_N) = \frac 1 \ZNA
\prod_{i <j} |x_i - x_j|^\beta I_N^\beta(\t, X_N) e^{-\frac N 2 \sum_{i=1}^N x_i^2}
dx_1 \ldots dx_N,$$
the largest eigenvalue $x_N^* = \max(x_1, \ldots, x_N)$ satisfies a
large deviation principle, in the scale $N$, with good rate function
$\KBT$ defined as follows:\\
$\bullet$ If $\t\ppq \sqrt{\frac\b 2},$ \begin{center}
\[
\KBT(x) = \left\{
\begin{array}{ll}
+\infty, & \textrm{ if } x< \sqrt{2\b}\\
\displaystyle \int_{\sqrt{2\b}}^x \sqrt{z^2-2\b} \,dz, & \textrm{ if } \sqrt{2\b} \ppq x
\ppq \t +\frac \b{2\t}, \\
M_\t^\b(x), &  \textrm{ if }  x \pgq  \t +\frac \b{2\t},

\end{array}
\right.
\]
with $\displaystyle  M_\t^\b(x) = \frac 1 2 \int_{\sqrt{2\b}}^x
\sqrt{z^2-2\b} \,dz
-\t x +\frac 1 4 x^2 + \frac \b 4 - \frac \b 4 \log \frac \b 2 + \frac
1 2 \t^2 + \frac{\b}{2} \log \t.$ \\
\end{center}
$\bullet$ If $\t\pgq \sqrt{\frac\b 2},$ \begin{center}
\[
\KBT(x) = \left\{
\begin{array}{ll}
+\infty, & \textrm{ if } x< \sqrt{2\b}\\
L_\t^\b(x), &  \textrm{ if }  x \pgq  \sqrt {2\b},
\end{array}
\right.
\]
with $\displaystyle  L_\t^\b(x) = \frac 1 2 \int_{ \t +\frac \b{2\t}}^x
\sqrt{z^2-2\b} dz - \t\left(x-\left( \t +\frac \b{2\t}\right)\right) 
+ \frac 1 4 \left(x^2-\left( \t +\frac \b{2\t}\right)^2\right).$
\end{center}
\end{theo}

The remaining of this section will be devoted to the proof of Theorem
\ref{main}.
A first step will be to prove the following LDP
\begin{prop} \label{ldp}
If we define 
\begin{equation} \label{PNA}
\PNA(dx_1, \ldots, dx_N) = \frac 1 {Z_N^\b} 
\prod_{i <j} |x_i - x_j|^\beta I_N^\beta(\t,X_N) e^{-\frac N 2 \sum_{i=1}^N x_i^2}
dx_1 \ldots dx_N,
\end{equation}
with $Z_N^\b$ the normalizing constant in the case $\t=0,$
 then 
under $\PNA$, the largest eigenvalue $x_N^*$ satisfies 
a large deviation principle, in the scale $N$, with good rate function
\[
F_\t^\b(x) := \left\{
\begin{array}{ll}
+\infty, & \textrm{ if } x< \sqrt{2\b} \\
- \frac \b 2 + \frac \b 2 \log \frac \b 2 - \Phi_\b(x, \s_\b) -I_{\s_\b}^\b(x,\t), & \textrm{ otherwise,}
\end{array}
\right.
\]
where $\s_\b$ denotes the semicircle law whose density on $\RR$ is
given by $\displaystyle \frac 1 {\b\pi} \mathbf{1}_{[-\sqrt{2\b},
    \sqrt{2\b}]} \sqrt{2\b - t^2} dt,$ \linebreak 
for  $\mu \in \PR$ and $x \in \RR$,
$$ \Phi_\b(x, \mu) = \b \int \log|x-y| d\mu(y) - \frac 1 2 x^2, $$
and $\displaystyle I_\mu^\b(x,\t) = \lim_{N \ra \infty} \frac 1 N \log
I_N^\b(\t, B_N),$ where $B_N
$
has limiting spectral measure $\mu$ and limiting largest eigenvalue $x.$ 
\end{prop}

\subsection{Proof of Proposition\ref{ldp} }
\mbox{}\\

The proof of Proposition \ref{ldp} will require the following
exponentiel tightness result :

 \begin{lem}
\label{tight}
For any $\t>0,$
there exists a function $f_\t: \RR^+ \ra \RR^+$ going to infinity at infinity
such that for all $N$
$$ \PNA\left( \max_{i=1 \ldots N} |x_i| \pgq M\right) \ppq e^{-N f_\t(M)} .$$
\end{lem}

\noindent
\textbf{Proof of Lemma \ref{tight}:}
It is more convenient to rewrite (\ref{PNA})
as
$$ \PNA(dx_1, \ldots, dx_N) = \frac {e^{\frac N 2 \t^2}}{Z_N^\b} 
\prod_{i <j} |x_i - x_j|^\b e^{-\frac N 2 \tr (X_N -A_N)^2}
dx_1 \ldots dx_N.$$
Now, a well known inequality (see for example Lemma 2.3 in \cite{Bai2})
gives that
$$ \tr (X_N -A_N)^2 \pgq \min_{\pi} \sum_{i=1}^N |x_k - a_{\pi(k)}|^2,$$ 
where the minimum is taken over all permutations $\pi$ of $\{1, \ldots, N\}$.
But all $a_k$'s are zero, except one of them, let's say $a_1$, which is equal to $\t$.
As the law of the $x_j$'s in invariant by permutations, we can assume that 
$\pi^{-1}_*(1)=1$, where $\pi_*$ is the permutation for which the minimum is reached.
Therefore
$$  \tr (X_N -A_N)^2 \pgq (x_1 -\t)^2 +  \sum_{i=2}^N  x_j^2.$$
We can now use the very same estimates as in Lemma 6.3 in \cite{BDG}
to get Lemma \ref {tight}. More precisely,  we can write
$$ |(x-\t) +\t - x_j|^\b e^{-\frac{x_j^2}{2}} \ppq e^{\frac{(x-\t)^2}{4}},$$
for $x$ large enough, so that, for $M$ large enough,
$$ \PNA\left( \max_{i=1 \ldots N} |x_i| \pgq M\right) \ppq N \PNA\left(  |x_1| \pgq M\right) 
\ppq \frac{Z_{N-1}^\b}{Z_N^\b} e^{- \frac 1 4 N (M-\t)^2 +\frac N 2 \t^2}.$$
From Selberg formula (cf for example proof of Proposition 3.1 in
\cite{BA-G}),
we can show that \newline
$\displaystyle \frac 1 N \log \frac{Z_{N-1}^\b}{Z_N^\b}
\xrightarrow[N \ra \infty]{} C.$
This concludes the proof of exponential tightness.\hfill $\Box$

\noindent 
We now go back to the proof of Proposition \ref{ldp}.\\

\noindent
$\bullet$ \textbf{$\FBT$ is a good rate function.}\\
From Theorem 6 in \cite{GM2}, it's not hard to check that $\KBT$
is continuous on $[\sqrt{2\b}, +\infty)$ and therefore lower
    semi-continuous on $\RR.$\\
Moreover, we can check that, for $x$ large enough, $I_{\s_\b}^\b(x,
\t) \ppq \t x,$ so that $\FBT(x) \sim_{+\infty} \frac 1 2 x^2,$
its level sets are therefore compact.  \\

\noindent
$\bullet$ For all $x < \sqrt{2\b}$,
\begin{equation}
\label{moinsdeux}
\lim_{N \ra \infty} \frac 1 N \log \PNA(x_N^* \ppq x) = - \infty.
\end{equation}
We know from Theorem 1.1 in \cite{BA-G} that the spectral measure
of $W_N$ satisfies a large deviation principle in the scale $N^2$
with a good rate function whose unique minimizer is the semicircle law
$\s_\b.$ 
We can check that adding a deterministic matrix of bounded rank
(uniformly in $N$) does not affect the spectral measure in this scale
so that the spectral measure of $X_N$ satisfies the same 
large deviation principle.\\
Therefore, if we let $x<\sqrt{2\b}$, $f \in \Ca_b(\RR)$ such that
$f(y)=0$ if $y \ppq x$ but $\int f d\s_\b >0$
and if we consider the closed set 
$F:= \{ \mu / \int f d\mu = 0 \}$, we have that
$$ \PNA(x_N^* \ppq x) \ppq \PNA\left(\frac 1 N \sum_{i=1}^N f(x_i)=0\right)
\ppq \PNA(\muN \in F),$$ 

where $\muN := \frac 1 N \sum_{i=1}^N \delta_{x_i}$ is the spectral measure
of $X_N$. As $\s_\b \notin F$,
$$ \limsup_{N \ra \infty} \frac{1}{N^2} \log \PNA(x_N^* \ppq x) <0,$$ 
what gives immediately (\ref{moinsdeux}).\\

\noindent
$\bullet$ Let now $x \pgq \sqrt{2\b}$ and $\d >0$. We want to show \textbf{the upper bound}.\\
Thanks to Lemma \ref{tight}, we can restrict ourselves 
to the event $\{ \max_{i=1 \ldots N} |x_i| \ppq M \},$
for an appropriate $M$.\\
One important remark is that, by invariance by permutation, we have,
for any real $x$,
$$ \PNA(x \ppq x_N^* \ppq x+\d, \max_{i=1 \ldots N} |x_i| \ppq M )
\ppq N \PNA(x \ppq x_1 \ppq x+ \d, \, x_1 \pgq \max_{i=2 \ldots N} x_i,
\max_{i=1 \ldots N} |x_i| \ppq M ).$$ 
We introduce now the following notations :
\begin{itemize}
\item $\displaystyle \pin := \frac{1}{N-1} \sum_{i=2}^N \delta_{x_i},$\\
\item $\PP_N^{N-1}$ is the measure on $\RR^{N-1}$ such that,
for each  Borelian set $E$, we have 
$$ \PP_N^{N-1}(\l \in E) = \PP_{N-1}^0\left( \sqrt{1-\frac 1 N} \l \in E\right).$$
\end{itemize} 
With these notations, we have 
\begin{multline*}
 B:= \PNA(x \ppq x_N^* \ppq x+\d, \max_{i=1 \ldots N} |x_i| \ppq M )\\
\ppq 
\int_x^{x+\d} dx_1 \int_{[-M,M]^{N-1}} e^{(N-1) \Phi_\b(x_1, \pin)}.
C_{N}^\b.I_N^\b(\t, X_N).
d\PP_{N-1}^0(x_2, \ldots, x_N),
\end{multline*}
where $\displaystyle C_{N}^\b :=
N \frac{Z_{N-1}^\b}{Z_{N}^\b}.\left(1-\frac 1 N\right)^{\b\frac{N(N-1)}{4}}.$\\

Let $0<\k<\frac 1 4$, we have
\begin{multline}
\label{maj}
B \ppq  C_{N}^\b.
\int_x^{x+\d} dx_1 \int_{\substack { \pin \in B(\s_\b, N^{-\k}),\\ \max_{i=2\ldots N} x_i \ppq x_1}}
 e^{(N-1) \Phi_\b(x_1, \pin)}
 I_N^\b(\t, X_N).
d\PP_N^{N-1}(x_2, \ldots, x_N)\\
+ (2M)^N e^{NM\t}  C_{N}^\b \PP_N^{N-1}(\pin \notin B(\s_\b, N^{-\k})),
\end{multline}
where $B(\s_\b, N^{-\k})$ is the ball of size $N^{-\k}$ centered at
$\s_\b$, 
for the Levy  distance. 

We first treat show that the second term  is exponentially negligible.
We have
$$ \PP_N^{N-1} (\pin \notin B(\s_\b, N^{-\k}))
\ppq \PP_N^{N-1}(\|F_{N-1} -F_\b\|_\infty \pgq N^{-\k}), $$
where $F_{N-1}$ and $F_\b$ are respectively the (cumulative)
distribution function of $\pin$ and $\s_\b$.\\
We know from the result of Bai 
in \cite{Bai3} that
$$ \|\E_N^{N-1} F_{N-1} - F_\b\| = O(N^{-\frac 1 4 }),$$
where $\EE_N^{N-1} $ is the expectation under $\PP_N^{N-1} $, so that
$$ \PP_N^{N-1}(\|F_{N-1} -F_\b\|_\infty \pgq N^{-\k})
\ppq \PP_N^{N-1}(\|F_{N-1} -\EE_N^{N-1}  F_{N-1}\|_\infty \pgq N^{-\k}).$$
But, by a result of concentration of \cite{GZ2} (see Theorem 1.1),
we have that there exists a constant $C>0$ such that for all $N \in \NN$,
$$ \PP_N^{N-1}(\|F_{N-1} -\EE_N^{N-1}  F_{N-1}\|_\infty \pgq N^{-\k})
\ppq e^{-C N^{2-2\k}},$$ 
so that
$$ \limsup_{N \ra \infty} \frac 1 N \log \PP_N^{N-1} (\pin \notin B(\s_\b, N^{-\k})) = -\infty.$$
We can now come back to the first term in (\ref{maj}).
The same computation as in the proof of Proposition 3.1 in
\cite{BA-G}, based on Selberg formula, gives that
$$ \frac 1 N \log C_N^\b \xrightarrow[N\ra \infty]{} -\frac \b 2 \log
\frac \b 2 + \frac \b 2.$$
 Applying Proposition \ref{cont} together with Theorem 6 of \cite{GM2}, we can conclude that
$$ \limsup_{N \ra \infty} 
\frac 1 N \log \PNA (x \ppq x_N^* \ppq x+\d, \max_{1 \ldots N}|x_i| \ppq M)
\ppq  -\frac \b 2 \log
\frac \b 2 + \frac \b 2
+ \sup_{z \in [x,x+\d]} [\Phi_\b(z,\s_\b) + I_{\s_\b}^\b(z,\t)].$$
As mentioned above,
$z \mapsto \Phi_\b(z, \s_\b) + I_{\s_\b}^\b(z, \t)$ is continuous on $[\sqrt{2\b}, +\infty)$
so that in particular
$$ \limsup_{\d \downarrow 0} \limsup_{N \ra \infty} 
\frac 1 N \log \PNA (x \ppq x_N^* \ppq x+\d, \max_{1 \ldots N}|x_i| \ppq M)
\ppq -\frac \b 2 \log
\frac \b 2 + \frac \b 2
+ \Phi_\b(x,\s_\b) + I_{\s_\b}^\b(x,\t).$$

\noindent
$\bullet$ We now conclude the proof of Proposition \ref{ldp} by showing the corresponding \textbf{lower bound}.
We proceed as in \cite{BDG}. Let $y>x>r> \sqrt{2\b}$. Then,
 \begin{eqnarray*}
\PNA(y \pgq x_N^* \pgq x) & \pgq & \PNA(x_1 \in [x,y], \max_{i=2\ldots N}|x_i| \ppq r)\\
& \pgq & C_N^\b  \exp\left( (N-1) 
\inf_{ \substack{
z \in [x,y] \\
\mu \in B_r(\s_\b, N^{-\k})
}} (\Phi_\b(z, \mu) + I_\mu(z, \t) - g_\k(x-y))\right)\\
& &.\PP_N^{N-1}(\nun \in B_r(\s_\b, N^{-\k})),
\end{eqnarray*}
where $B_r(\s_\b, N^{-\k}) = B(\s_\b, N^{-\k}) \cap \Pa([-r,r])$ and $g_\k$
going to zero at zero
by virtue of Proposition \ref{cont}.\\
We proceed as for the upper bound to show that
$\PP_N^{N-1}(\pin \in B_r(\s_\b, N^{-\k}))$ is going to $1$.  
Knowing the asymptotics of $C_N^\b$, we get
$$ \liminf_{N \ra \infty} \PNA(y \pgq x_N^* \pgq x) 
\pgq  -\frac \b 2 \log
\frac \b 2 + \frac \b 2 + \inf_{
z \in [x,y]} (\Phi_\b(z, \s_\b) + I_{\s_\b}^\b(z, \t)-g_\k(x-y)). $$
We let now $y$ decrease to $x$. $\Phi_\b(.,\s_\b)$ and $I_{\s_\b}^\b(.,\t)$
are continuous on $(\sqrt{2\b}, +\infty[$ (we are outside the support of
$\s_\b$) so that we have the required lower bound
$$ \liminf_{y \ra x} \liminf_{N \ra \infty} \PNA(x_N^* \pgq x) 
\pgq -\frac \b 2 \log
\frac \b 2 + \frac \b 2 +  \Phi_\b(x, \s_\b) + I_{\s_\b}^\b(x, \t).$$ 
 This concludes the proof of Proposition  \ref{ldp}. \hfill $\Box$\\

\subsection{Proof of Theorem \ref{main}}
\mbox{}\\

\noindent
A direct consequence of Proposition  \ref{ldp} is that, for
$$ \ZNA = \int \ldots \int \prod_{i <j} |x_i - x_j|^\beta I_N^\beta(\t,X_N) e^{-\frac N 2 \sum_{i=1}^N x_i^2}
dx_1 \ldots dx_N,$$
we have that  $\displaystyle \lim_{N\ra \infty}\frac 1 N \log \frac \ZNA {Z_N^\b} =
\inf_{x \in \RR}
\FBT(x)$ so that $\KBT(x) = \FBT(x) - \inf_{x \in \RR}
\FBT(x).$\\
To conclude the proof of Theorem \ref{main}, we have to study the
function $\FBT.$ We recall that, for $x \pgq \sqrt{2\b},$ we have
$$ \FBT(x) = - \frac \b 2 + \frac \b 2 \log \frac \b 2 - \b \int 
\log|x-y| d\s_\b(y) +\frac 1 2 x^2 -I_{\s_\b}^\b(x,\t),$$
where we know from Theorem 6 in \cite{GM2} that
 $$I_{\mu}^\b(x,\t) = \t v(x, \t) -\frac \b 2 \int \log\left(1+\frac{2\t}\b v(x, \t) - \frac{2\t}\b \l\right) d\mu(\l),$$
\[
\textrm{ with } \quad v(x, \t) := \left\{
\begin{array}{ll}
\R\left(\frac{2\t}\b\right), & \textrm{ if } \H(x) \pgq \frac{2\t}\b, \\
x - \frac{\b}{2\t}, & \textrm{ otherwise,}
\end{array}
\right. .
\]

Relying for example on the proof of Lemma 2.7 in \cite{BA-G}, we have
that
\begin{equation}\label{H2}
 - \frac \b 2 + \frac \b 2 \log \frac \b 2 - \b \int
\log|x-y| d\s_\b(y) +\frac 1 2 x^2 = \int_{\sqrt{2\b}}^x \sqrt{z^2-2\b
} dz,
\end{equation}
 so that 
\begin{equation}\label{H3}
 \FBT(x) =  \int_{\sqrt{2\b}}^x \sqrt{z^2-2\b
} dz - I_{\s_\b}^\b(x,\t).
\end{equation}
From  Lemma 2.7 in \cite{BA-G}, we also get that
\begin{equation}\label{H}
 H_{\s_\b}(x) = \frac 1 \b (x-\sqrt{x^2-2\b }) \textrm{ for }
x>\sqrt{2\b},
\end{equation}
from what we deduce 
that $\displaystyle \lim_{x \ra \sqrt{2\b}}
H_{\s_\b}(x) = \frac {\sqrt{2\b}} \b = \sqrt{\frac 2 \b}$
and $H_{\s_\b}$ is decreasing.\\

\noindent
$\bullet$ If $\t \pgq \sqrt{\frac \b 2},$ then for all $x>\sqrt{2\b}$,
$\displaystyle  H_{\s_\b}(x) \pgq \frac{2\t}\b$ and
$$\displaystyle I_{\s_\b}^\b(x,\t) = \t x - \frac \b 2  - \frac \b 2 \log \left( \frac{2\t}\b\right)
- \frac \b 2 \int \log |x-y| d\s_\b(y) :=S_\t(x),$$
so that from \eqref{H2} and \eqref{H3}, we get that 
$$\displaystyle \FBT(x) = \frac  1 2  \int_{\sqrt{2\b}}^x \sqrt{z^2-2\b
} dz  -\t x + \frac 1 4 x^2 + \frac \b 4  - \frac \b 2 \log \frac \b 2
+ \frac \b 2 \log \t.$$\\
Differentiating this function on $(\sqrt{2\b}, +\infty),$
we see that it is decreasing on $(\sqrt{2\b}, \t +\frac\b{2\t})$ and
  then increasing so that its infimum is reached at $\t
  +\frac\b{2\t}.$
This gives immediately  in this case that 
$\KBT(x)=\FBT(x) - \inf_x \FBT(x) = L_\t^\b(x),$ as defined in Theorem
\ref{main}.\\

\noindent
$\bullet$  If   $\t \ppq \sqrt{\frac  \b 2},$ then we can check that 
on   $\left[\sqrt{2\b}, \t +\frac\b{2\t}\right],$ we have 
$\displaystyle  H_{\s_\b}(x) \pgq \frac{2\t}\b$ and \newline
$\displaystyle I_{\s_\b}^\b(x,\t) = \frac \b 2 \int_0^{\frac{2\t}\b}
R_{\s_\b}(u)du.$
Moreover, from  \eqref{H}, we get that the inverse of $H_{\s_\b}$
is given by $G_{\s_\b}(x) = \frac \b 2 x +\frac 1 x$ so that
$R_{\s_\b}(x) =  \frac \b 2 x $ and 
$ I_{\s_\b}^\b(x,\t) = \frac 1 2 \t^2.$\\
In this case, $\displaystyle  \FBT(x) =  \int_{\sqrt{2\b}}^x \sqrt{z^2-2\b
} dz - \frac 1 2 \t^2, $ which is increasing.\\
For $x> \t +\frac\b{2\t}, \displaystyle  H_{\s_\b}(x) \ppq \frac{2\t}\b,$ 
so that $ I_{\s_\b}^\b(x,\t) = S_\t(x)$ as above.\\
Therefore, $\FBT$ is increasing on  $\left[\sqrt{2\b}, \t +\frac\b{2\t}\right]$
and on  $\left[\t +\frac\b{2\t}, +\infty\right)$ and is continuous
so that its infimum is reached at $\sqrt{2\b}$
and is equal to $-\frac 1 2 \t^2.$\\
Therefore, $\displaystyle \KBT(x) =  \int_{\sqrt{2\b}}^x \sqrt{z^2-2\b
} dz$ on  $[\sqrt{2\b}, \t +\frac\b{2\t}]$
and coincides with $M_\t^\b(x)$ on  $\left[\t +\frac\b{2\t}, +\infty\right).$
This concludes the proof of Theorem \ref{main}. \hfill $\Box$

\section{Proof of Corollary \ref{seppgv}}
\label{sec:sep}

In the proof of Theorem \ref{main} above, we proved that $ \KBT$
is increasing on  $[\sqrt{2\b},+\infty)$ if $\t \ppq \sqrt{\frac \b
      2},$
so that in this case its infimum is reached at $\sqrt{2\b}.$\\
We also saw that, differentiating $\FBT$ on   $(\sqrt{2\b},+\infty),$
we got that when $\t \pgq \sqrt{\frac \b
      2},$ it reaches its minimum at $ \t +\frac\b{2\t}.$
This is enough to conclude. \hfill $\Box$
\mbox{}\\

\noindent
{\bf Acknowledgments :} I would like to thank Sandrine P\'ech\'e for many
  fruitful discussions during the preparation of the paper. We also would like to thank G\'erard Ben Arous for pointing out a typo in the expression of the rate function in Theorem \ref{main}. We are very indebted to Benjamin McKenna, who found a mistake in the proof of Lemma 2.3 in the published version of the paper and proposed a new proof (through Lemma 2.5). \\

\bibliographystyle{acm}
\bibliography{rank1}

\begin{thebibliography}{10}

\bibitem{Bai3}
{\sc Bai, Z.~D.}
\newblock Convergence rate of expected spectral distributions of large random
  matrices. {I}. {W}igner matrices.
\newblock {\em Ann. Probab. 21}, 2 (1993), 625--648.

\bibitem{Bai2}
{\sc Bai, Z.~D.}
\newblock Methodologies in spectral analysis of large-dimensional random
  matrices, a review.
\newblock {\em Statist. Sinica 9}, 3 (1999), 611--677.
\newblock With comments by G. J.\ Rodgers and Jack W.\ Silverstein; and a
  rejoinder by the author.

\bibitem{BDG}
{\sc Ben~Arous, G., Dembo, A., and Guionnet, A.}
\newblock Aging of spherical spin glasses.
\newblock {\em Probab. Theory Related Fields 120}, 1 (2001), 1--67.

\bibitem{BA-G}
{\sc Ben~Arous, G., and Guionnet, A.}
\newblock Large deviations for {W}igner's law and {V}oiculescu's
  non-commutative entropy.
\newblock {\em Probab. Theory Related Fields 108}, 4 (1997), 517--542.

\bibitem{BH}
{\sc Br{\'e}zin, E., and Hikami, S.}
\newblock Correlations of nearby levels induced by a random potential.
\newblock {\em Nuclear Phys. B 479}, 3 (1996), 697--706.

\bibitem{DKM}
{\sc Deift, P., Kriecherbauer, T., McLaughlin, K. T.-R., Venakides, S., and
  Zhou, X.}
\newblock Uniform asymptotics for polynomials orthogonal with respect to
  varying exponential weights and applications to universality questions in
  random matrix theory.
\newblock {\em Comm. Pure Appl. Math. 52}, 11 (1999), 1335--1425.

\bibitem{Dei}
{\sc Deift, P.~A.}
\newblock {\em Orthogonal polynomials and random matrices: a
  {R}iemann-{H}ilbert approach}, vol.~3 of {\em Courant Lecture Notes in
  Mathematics}.
\newblock New York University Courant Institute of Mathematical Sciences, New
  York, 1999.

\bibitem{FP}
{\sc F{\'e}ral, D., and P{\'e}ch{\'e}, S.}
\newblock The largest eigenvalue of rank one deformation of large wigner
  matrices.
\newblock Preprint, math.PR/0605624, 2006.

\bibitem{GM2}
{\sc Guionnet, A., and Ma{\"{\i}}da, M.}
\newblock A {F}ourier view on the {$R$}-transform and related asymptotics of
  spherical integrals.
\newblock {\em J. Funct. Anal. 222}, 2 (2005), 435--490.

\bibitem{GZ}
{\sc Guionnet, A., and Zeitouni, O.}
\newblock Large deviations asymptotics for spherical integrals.
\newblock {\em J. Funct. Anal. 188}, 2 (2002), 461--515.

\bibitem{GZ2}
{\sc Guionnet, A., and Zeitouni, O.}
\newblock Addendum to large deviations asymptotics for spherical integrals.
\newblock {\em J. Funct. Anal. To appear\/} (2004).

\bibitem{HR}
{\sc Hoyle, D., and Rattray, M.}
\newblock Limiting form of the sample covariance eigenspectrum in pca and
  kernel pca.
\newblock In {\em Proceedings of Neural Information Processing Systems\/}
  (2003).

\bibitem{Joh}
{\sc Johansson, K.}
\newblock Universality of the local spacing distribution in certain ensembles
  of {H}ermitian {W}igner matrices.
\newblock {\em Comm. Math. Phys. 215}, 3 (2001), 683--705.

\bibitem{LCPB}
{\sc Laloux, L., Cizeau, P., Potters, M., and Bouchaud, J.}
\newblock Random matrix theory and financial correlations.
\newblock {\em Intern. J. Theor. Appl. Finance 3}, 3 (2000), 391--397.

\bibitem{Meh}
{\sc Mehta, M.~L.}
\newblock {\em Random matrices}, second~ed.
\newblock Academic Press Inc., Boston, MA, 1991.

\bibitem{SP}
{\sc P{\'e}ch{\'e}, S.}
\newblock The largest eigenvalue of small rank perturbations of hermitian
  random matrices.
\newblock Preprint, math.PR/0411487, 2004.

\bibitem{TW}
{\sc Tracy, C.~A., and Widom, H.}
\newblock Level-spacing distributions and the {A}iry kernel.
\newblock {\em Comm. Math. Phys. 159}, 1 (1994), 151--174.

\bibitem{Wig}
{\sc Wigner, E.~P.}
\newblock On the distribution of the roots of certain symmetric matrices.
\newblock {\em Ann. of Math. (2) 67\/} (1958), 325--327.

\end{thebibliography}
\end{document}